\documentclass[a4paper,10pt,fleqn]{article}
\usepackage{graphics,graphicx}
\usepackage{natbib}
\usepackage{amsmath}
\usepackage{amssymb}

\chardef\at=`\@
\DeclareRobustCommand{\qed}{%
  \ifmmode 
  \else \leavevmode\unskip\penalty9999 \hbox{}\nobreak\hfill
  \fi
  \quad\hbox{\qedsymbol}}
\newcommand{\mathbold}[1]{\mbox{\boldmath $#1$}}

\newcommand{\pr}{\mathbold P}
\newcommand{\ep}{\mathbb P}

\newcommand{\ex}{\mathbold E}

\newcommand{\argmin}{\operatorname{argmin}}

\begin{document}
\begin{center}
{\large Functional Choice and Non-significance Regions in Regression\\}
\quad\\
Laurie Davies\\
Fakul\"at Mathematik\\
Universit\"at Duisburg-Essen\\
\texttt{laurie.davies@uni-due.de}
\end{center}
\quad\\
\begin{abstract}
Given data $y$ and $k$ covariates $x$ the problem is to
decide which covariates to include when approximating $y$ by a linear
function of the covariates. The decision is based on replacing subsets
of the covariates by i.i.d. normal random variables and comparing the
error with that obtained by retaining the subsets. If the two errors
are not significantly different for a particular subset it is
concluded that the covariates in this subset are no better than random
noise and they are not included in the linear approximation to $y$. 
\end{abstract}

\section{Introduction} \label{sec:intro}
\subsection{Notation}
Consider $n$ measurements ${\mathbold y}_n=(y_1,\ldots,y_n)^t$  of a
variable $y$ and for each $y_i$ concomitant measurements of $k$
covariables $x_j, j=1,\ldots,k,$ given by  ${\mathbold
  x}_{i\cdot}=(x_{i1},\ldots,x_{ik})^t$ forming an $n\times k$ matrix
${\mathbold x}_n$ with $j$th column ${\mathbold x}_{\cdot j}$. A
subset of the covariates will be denoted by a row
vector $e=(e_1,\ldots,e_k)$ with $e_j\in \{0,1\}$ whereby $e_j=1$
means that the $j$th covariate is included.  A model $e$ will be
encoded as $\sum_{j=1}^ke_j2^{j-1}$. The subset consisting of
all covariates will be denoted by $e_f$. Given an $e$ with
$\sum_{j=1}^k e_j=k(e)$ the $n\times k(e)$ matrix with columns
corresponding to those covariates with $e_j=1$ will be 
denoted by ${\mathbold x}_n(e)$ with ${\mathbold x}_{i\cdot}(e)$,
${\mathbold x}_{\cdot j}(e)$ and ${\mathbold x}(e)$ having the
corresponding interpretations. The empirical measure of the data will
be denoted by $\ep_n$
\begin{equation} \label{equ:emp_measure}
\ep_n=\ep_n(({\mathbold y}_n,{\mathbold
  x}_n))=\frac{1}{n}\sum_{i=1}^n\delta_{y_i,{\mathbold x}_{i\cdot}}
\end{equation}
with the corresponding definition of $\ep_{n,e}$ for any subset $e$. The
$L_1$ and $L_2$ norms will be denoted by $\Vert\,\Vert_1$ and
$\Vert\,\Vert_2$ respectively.

\subsection{The problem}
 The problem to decide which
if any of the covariates $x_j$  influence the value of $y$. There
are many proposals for doing this. Some such  as AIC
(\cite{AKAIK73,AKAIK74,AKAIK81}) or BIC (\cite{SCH78}) require an
explicit model such as 
\begin{equation} \label{equ:linreg1_mod}
Y={\mathbold x}^t{\mathbold \beta}+\varepsilon
\end{equation}
where ${\mathbold x}^t=(x_1,\ldots,x_k)$, ${\mathbold
  \beta}=(\beta_1,\ldots,\beta_k)^t$ and the errors
$\varepsilon$ are  random variables with an explicit
distribution. Others such as Lasso (\cite{TIB96})
\begin{equation} \label{equ:linreg1}
\argmin_{\mathbold \beta}\left\{\sum_{i=1}^n (y_i-{\mathbold
    x}_{i\cdot}{\mathbold \beta})^2+\lambda \sum_{j=1}^k \vert
  \beta_j\vert \right\}
\end{equation}
may or may not require an explicit model to determine the choice of
the smoothing parameter $\lambda$. 

The following is based on a simple idea.   Let $s^2$ denotes the least
sum of squares based on all covariates and for a given subset let
$S_e^2$ denote the least sum of squares when all the covariates with
$e_j=0$ are  replaced by i.i.d. $N(0,1)$ random variables. If $s^2$ is
not significantly small than $S_e^2$ the conclusion is that the
omitted covariates are no better than random noise. To define
`significantly' the process is repeated a large number of times. For a
given $\alpha$, $\alpha=0.95$ for example, $s^2$ is significantly
smaller the $S_e^2$ if in at least $100\alpha\%$ the simulations
$s^2\le S_e^2$. The P-value $p_e$ is the proportion of simulations for
which  $S_e^2 < s^2$, so that the excluded covariates are
significantly better than random noise if $p_e \le 1-\alpha$. A small
values of $p_e$ indicates that at least some of the omitted covariates
are relevant. A large value of $p_e$ indicates that in toto the
omitted covariates are no better than random noise.

The method is not restricted to least squares regression. It can be
equally well applied to $L_1$ regression or more generally to any
measure of discrepancy  $d({\mathbold y}_n,{\mathbold x}_n)$. 

As an example consider the stack loss data of \cite{BRO60}. It is
one of the data sets provided by \cite{R13} and is used in
\cite{KOEN10}. There are 21 observations with one dependent variable 
`Stack.Loss' and the three covariates `Air.Flow', `Water.Temp' and
`Acid.Conc' labelled from one to three. In the following the intercept
will always be included. There are eight possible models. The results of an $L_1$
regression for the stack loss data are given in
Table~\ref{tab:stackloss_L_1}. The total computing time was 102
seconds using \cite{KOEN10}. The only subset with a large $P$-value is
the subset encoded as 3 which corresponds to $e=(1,1,0)$. 
\begin{table}[h]
\begin{tabular}{ccccccccc}
subset&0&1&2&3&4&5&6&7\\
$P$-value&0.000&0.015&0.000&0.231&0.000&0.007&0.000&1.000\\
\quad\\
\end{tabular}
\caption{Encoded subsets and $P$-values for the stack loss
  data based on 5000 simulations.\label{tab:stackloss_L_1}}
\end{table}

\subsection{Non-significance regions} \label{sec:n-s-r_example}
Given a subset $e$ of covariates the best linear approximation to the
variable ${\mathbold y}_n$ in the $L_1$ norm is
\begin{equation} \label{equ:best_l_1}
{\mathbold x}_n(e){\mathbold \beta}_{1,n}(e)
\end{equation}
where
\begin{equation} \label{equ:best_l_1_1}
{\mathbold \beta}_{1,n}(e)= \argmin_{{\mathbold \beta}(e)}\Vert
{\mathbold y}_n-{\mathbold x}_n(e){\mathbold \beta}(e)\Vert_1\,.
\end{equation}
 A single value is not sufficient to answer many questions of
interest which require a range of plausible values.  In frequentist
statistics such a range is provided by a confidence region. This
option is not available in the present context as a confidence region
assumes that there is a `true' value to be covered.  The confidence
region will be replaced by a non-significance region whose
construction will be illustrated for the median. 

Given data ${\mathbold y}_n$ the median minimizes $s_1({\mathbold y}_n)=\sum_{i=1}^n\vert
 y_i-\text{med}({\mathbold y}_n)\vert $. For any other value $m\ne
 \text{med}({\mathbold y}_n)$ 
\[\sum_{i=1}^n\vert y_i-\text{med}({\mathbold y}_n)\vert <
\sum_{i=1}^n\vert y_i-m\vert \]
A value $m$ will be considered as not being significantly different from
the median $\text{med}({\mathbold y}_n)$ if the difference
\[\sum_{i=1}^n\vert y_i-m\vert-\sum_{i=1}^n\vert
y_i-\text{med}({\mathbold y}_n)\vert\]
is of the order attainable by a random perturbation of the
$y$-values. More precisely if
\begin{equation} \label{equ:def_nonsig_med}
\pr\left(\inf_b\sum_{i=1}^n\vert y_i+bZ_i-m\vert < s_1({\mathbold
    y}_n)\right)\ge 1-\alpha.
\end{equation}
The set of values $m$ which satisfy (\ref{equ:def_nonsig_med}) can be
determined by simulations. For the ${\mathbold y}_n$ of the stack loss
data the 0.95-non-significance region is $[11.86,18.71]$ which can be
compared with the 0.95-confidence region $[11,18]$ based on the order
statistics. For any $m$ the P-value $p(m)$ is defined as
\begin{equation} \label{equ:def_nonsig_med-Pval}
p(m)=\pr\left(\inf_b\sum_{i=1}^n\vert y_i+bZ_i-m\vert < s_1({\mathbold
    y}_n)\right).
\end{equation}

\section{Choice of functional}
The procedure expounded in the previous section makes no use of a
model of the form (\ref{equ:linreg1_mod}). It solely based on the
approximation of ${\mathbold y}_n$ by a linear combination of the
covariates as measured in the $L_1$ and $L_2$ norms. There is no
mention of an error term. It therefore makes little sense to describe
the procedure as one of model sense. It makes more sense to interpret
it as one of functional choice. There does not seem to be any
immediate connection with 'wrong model' approaches as in
\cite{BERetal13} and \cite{LINLIU09}.

For a given subset $e$ the $L_1$ function
$T_{1,e}$ is defined by
\begin{eqnarray} \label{equ:functional_L_1}
T_{1,e}(\ep_n) &=& \argmin_{{\mathbold \beta}(e)}
\int\vert y-{\mathbold x}(e)^t{\mathbold \beta}(e)\vert \,
d\ep_n(y,{\mathbold x}(e))\nonumber\\
&=&\argmin_{{\mathbold \beta}(e)}\Vert
{\mathbold y}_n-{\mathbold x}_n(e){\mathbold \beta}(e)\Vert_1
\end{eqnarray}
with the corresponding definition of the $L_2$ functional
\begin{eqnarray} \label{equ:functional_L_2}
T_{2,e}(\ep_n) &=& \argmin_{{\mathbold \beta}(e)}
\int (y-{\mathbold x}(e)^t{\mathbold \beta}(e))^2 \,
d\ep_n(y,{\mathbold x}(e))\nonumber\\
&=&\argmin_{{\mathbold \beta}(e)}\Vert
{\mathbold y}_n-{\mathbold x}_n(e){\mathbold \beta}(e)\Vert_2.
\end{eqnarray}
More generally an $M$-functional $T_{\rho,e}$ can be defined as
\begin{eqnarray} \label{equ:functioal_rho_0}
T_{\rho,e}(\ep_n) &=& \argmin_{{\mathbold \beta}(e)}
\int\rho\left(\frac{y-{\mathbold x}(e)^t{\mathbold
      \beta}(e)}{\sigma_n}\right) \, dP(y,{\mathbold x}(e))\nonumber\\
&=&\argmin_{{\mathbold \beta}(e)}
\frac{1}{n}\sum_{i=1}^n\rho\left(\frac{y_i-{\mathbold
      x}_{i\cdot}(e)^t{\mathbold \beta(e)}}{\sigma_n}\right). 
\end{eqnarray}
The function $\rho$ is taken to be convex with a bounded first
derivative. This is the case for the  default choice in this paper
namely the Huber $\rho$-function defined by  
\begin{equation} \label{equ:huber_rho}
\rho_c(u)=\left\{\begin{array}{l@{\quad:\quad}l}
\frac{1}{2}u^2& \vert u\vert \le c\\
    c\vert u\vert-\frac{1}{2}c^2 &\vert u\vert > c\\
\end{array}\right.
\end{equation}
where $c$ is a tuning constant. The functional can be calculated using
the iterative scheme described in Chapter 7.8.2 of \cite{HUBRON09}.  

For reasons of equivariance (\ref{equ:functioal_rho_0}) contains a
scale parameter $\sigma_n$  which may be external or part of the
definition of $T_{\rho}$ (see Chapter 7.8 of \cite{HUBRON09}). The
default choice in this paper is the Median Absolute Deviation of the
residuals from an $L_1$ fit: 
\begin{equation} \label{equ:default_sig_2}
\sigma_n=\text{mad}({\mathbold y}_n-{\mathbold x}_n{\mathbold
  \beta}_{1,n}(e_f)).
\end{equation}
One use of $M$-functionals is to protect against outlying
$y$-values. The choice  (\ref{equ:default_sig_2}) preserves this
property.

\subsection{$L_1$ regression}\label{sec:rq_modchc}
The best linear fit based on all covariates is determined by
\begin{equation} \label{equ:rq_functional_1}
T_{1,e_{\text f}}(\ep_n)={\mathbold \beta}_{1,n}(e_{\text f})=
\argmin_{{\mathbold \beta}}\frac{1}{n}\sum_{i=1}^n\vert
y_i-{\mathbold x}_{i\cdot}^t{\mathbold \beta}\vert=\argmin_{{\mathbold
    \beta}}\Vert {\mathbold y}_n-{\mathbold x}_n{\mathbold \beta}\Vert_1 
\end{equation}
with mean sum of absolute deviations 
\begin{equation} \label{equ:rq_functional_3}
s_{1,n}(e_{\text f})=\frac{1}{n}\sum_{i=1}^n\vert y_i-
{\mathbold x}_{i\cdot}^t{\mathbold \beta}_{1,n}(e_{\text f})\vert =\Vert {\mathbold y}_n-
{\mathbold  x}_n{\mathbold \beta}_{1,n}(e_{\text f})\Vert_1 \,.
\end{equation}

 Let ${\mathbold Z}_n$ be a $n\times k$ matrix with elements $Z_{ij}$
which are i.i.d. $N(0,1)$. Given $e$ replace the covariates with
$e_j=0$ by the $Z_{ij}$, that is, put $W_{i,j}(e)=x_{i,j}$ if
$e_j=1$ and $W_{ij}(e)=Z_{ij}$ if $e_j=0$. Denote the relevant
matrices by ${\mathbold W}_n(e)$ and ${\mathbold
  Z}_n(e^c)$ and the empirical measure by ${\tilde \ep}_{n,e}$. The
best linear fit based on these covariates is determined  
\begin{equation}  \label{equ L_1_all}
T_{1,e_{\text f}}({\tilde \ep}_{n,e})={\tilde
  {\mathbold \beta}}_{1,n}(e)=\argmin_{{\mathbold \beta}} \Vert {\mathbold y}_n-
{\mathbold W}_n(e){\mathbold \beta}\Vert_1
\end{equation}
with mean sum of absolute deviations
\begin{equation} \label{equ:rq_functional_2}
S_{1,n}(e)=\Vert {\mathbold y}_n-{\mathbold W}_n(e) {\tilde
  {\mathbold \beta}}_{1,n}(e)\Vert_1
\end{equation}
The quantity $S_{1,n}(e)$ is a random variable. The $P$-value
$p_n(e)$ is defined by
\begin{equation} \label{equ:rq_functional_21}
p_n(e)={\mathbold P}(S_{1,n}(e)\le s_{1,n}(e_{\text f}))\,.
\end{equation}

There is no explicit expression for the $P$-values in the case of
$L_1$ regression. They must be calculated using simulations as in
Table~\ref{tab:stackloss_L_1}. This results in four of the
$P$-values being zero and so no comparison between them. A comparison 
can be obtained as follows. Simulate the distribution of 
\begin{equation} \label{equ:s_S_1}
s_{1,n}(e)- S_{1,n}(e)
\end{equation}
where
\begin{equation} \label{equ:s(e)}
s_{1,n}(e)=\frac{1}{n}\sum_{i=1}^n\vert y_i-
{\mathbold x}_{i\cdot}^t{\mathbold \beta}_{1,n}(e)\vert =\Vert {\mathbold y}_n-
{\mathbold  x}_n{\mathbold \beta}_{1,n}(e)\Vert_1 
\end{equation}
and then approximate it by a $\Gamma$-distribution with the shape and
scale parameters $sh(e)$ and $sc(e)$  estimated from the simulations
as $\hat{sh}(e)$ and $\hat{ sc}(e)$ respectively. The resulting estimated
$P$-values are given by
\begin{equation} \label{equ:p_gamma_approx_1}
{\hat p}_n(e)=1-\text{pgamma}(s_{1,n}(e)-s_{1,n}(e_{\text f}),
\hat{sh}(e),\hat{sc}(e))\,.
\end{equation}
The results for the stack loss data are given in
Table~\ref{tab:stack_loss_gamma} and may be compared with the
$P$-values of Table~\ref{tab:stackloss_L_1}.
\begin{table}[h]
\begin{tabular}{ccccccccc}
functional&0&1&2&3&4&5&6&7\\
$P$-value&1.93e-7&1.41e-2&4.90e-4&2.32e-1&5.02e-9&7.43e-3&2.57e-4&1.00\\
\quad\\
\end{tabular}
\caption{Encoded $L_1$-functionals and $P$-values for the stack loss
  data based on 1000 simulations using the $\Gamma$-approximation
  (\ref{equ:p_gamma_approx_1}),\label{tab:stack_loss_gamma}} 
\end{table}

Small $P$-values indicate that covariables have been omitted which have
a significant effect on the dependent variable. This excludes the
functionals encoded as 0, 1, 2, 4, 5, 6 although the functional
encoded as 1 could possible be retained. The functional 3 with
$P$-value 0.232 omits the covariable Acid.Conc.  As the functional 7
differs from 3 only through the inclusion of Acid.Conc the conclusion
is that it contains a covariate which is little  better than random
noise. Thus an analysis of the $P$-values leads to the choice of  the
functional 3. The interpretation of $P$-values and the choice of
functional will be considered in greater detail in
Sections~\ref{sec:eval_p} and \ref{sec:choose_func} respectively.

The second running example is the low birth weight data of
\cite{HOSLEM89} with $n=189$ and $k=9$. The dependent variable is the
weight of the child at birth. The nine covariates range from the
weight and age of the mother to hypertension and indicators of
race. There are in all 512 different functionals.  In the context of
model choice it is considered in \cite{HJOCLA03B}.

For this data set the computing time using 1000 simulations is about 50
minutes. This can be reduced by a factor of about ten by approximating
the modulus function $\vert x\vert$ by the Huber $\rho$-function
(\ref{equ:huber_rho}) with a small value of the tuning constant $c$,
for example $c=0.01$ (see Section~\ref{sec:M_modchc}). Care must be
taken in interpreting the decrease in computing time as the
$L_1$-functional was calculated using package \cite{KOEN10} whereas
the program for the $M$-functional was  written entirely in Fortran using
the algorithm given in Chapter 7.8  of \cite{HUBRON09} (see also
\cite{DUTT77a} and \cite{DUTT77b}) and  the pseudo-random number
generator {\it ran2} (see  \cite{PRTEVEFL03}). A pure Fortran program
for the $L_1$-functional  may be much faster (see
\cite{KoenkerPortnoy1997}).

\subsection{$M$-regression functionals}\label{sec:M_modchc}
The $M$-functionals can be treated in the same manner as the $L_1$
functional but with the added advantage that for large values of the
tuning constant $c$ in (\ref{equ:huber_rho}) there exist asymptotic
approximations for the $P$-values. On writing
\begin{eqnarray}
T_{\rho,e_{\text f}}(\ep_n) &=&{\mathbold \beta}_{\rho}(e_{\text f})=\argmin_{{\mathbold \beta}}
\frac{1}{n}\sum_{i=1}^n\rho\left(\frac{y_i-{\mathbold
      x}_{i\cdot}^t{\mathbold \beta}}{\sigma_n}\right)\\ 
s_{\rho,n}(e_{\text{f}})&=&\frac{1}{n}\sum_{i=1}^n\rho
\left(\frac{y_i-{\mathbold x}_{i\cdot}^t{\mathbold \beta}_{\rho}(e_{\text f})}{\sigma_n}\right)\\
T_{\rho,e}(\ep_n) &=&{\mathbold \beta}_{\rho}(e)=\argmin_{{\mathbold
    \beta}}\frac{1}{n}\sum_{i=1}^n 
\rho\left(\frac{y_i-{\mathbold x}_{i\cdot}(e)^t{\mathbold \beta}}{\sigma_n}\right)\\
s_{\rho}(e)&=&\frac{1}{n}\sum_{i=1}^n\rho\left(\frac{y_i-{\mathbold
      x}_{i\cdot}(e)^t{\mathbold \beta}_{\rho}(e)} 
{\sigma_n}\right)\\
T_{\rho,e}({\tilde \ep}_n) &=&{\tilde {\mathbold \beta}}_{\rho}(e)=\argmin_{{\mathbold \beta}(e)}
\frac{1}{n}\sum_{i=1}^n\rho\left(\frac{y_i-{\mathbold
      W}_{i\cdot}(e)^t{\mathbold \beta}(e)}{\sigma_n}\right)\\ 
S_{\rho,n}(e)&=&\frac{1}{n}\sum_{i=1}^n\rho\left(\frac{y_i-{\mathbold W}_{i\cdot}(e)^t
{\tilde {\mathbold \beta}}_{\rho}(e)}{\sigma_n}\right)
\end{eqnarray}
a second order Taylor expansion  gives
\begin{equation} \label{equ:taylor_approx_1}
S_{\rho,n}(e)\approx s_{\rho,n}(e)-\frac{1}{2}\frac{ \left(\frac{1}{n}
\sum_{i=1}^n\rho^{(1)}\left(\frac{r_i(e)}{\sigma_n}\right)^2\right)
\chi^2_{k-k(e)}}{\frac{1}{n}\sum_{i=1}^n\rho^{(2)}
\left(\frac{r_i(e)}{\sigma_n}\right)}
\end{equation}
where $\rho^{(1)}$ and $\rho^{(2)}$ are first and second derivatives
of $\rho$ respectively, $r_i(e)=y_i-{\mathbold x}_{i\cdot}(e)^t{\mathbold \beta}(e)$ and
$\chi^2_{k-k(e)}$ is a chi-squared random variable with $k-k(e)$
degrees of freedom. The inequality $S_{\rho,n}(e) \le s_{\rho,n}(e_{\text
  f})$ is asymptotically equivalent to
\begin{equation}
 \frac{\left(\frac{2}{n}\sum_{i=1}^n\rho^{(2)}\left(\frac{r_i(e)}{\sigma_n}
\right)\right)(s_{\rho,n}(e)-s_{\rho,n}(e_{\text
f}))}{\frac{1}{n}\sum_{i=1}^n
\rho^{(1)}\left(\frac{r_i(e)}{\sigma_n}\right)^2}\le \chi^2_{k-k(e)}
 \end{equation}
with asymptotic $P$-value
\begin{equation} \label{equ:asymp_p_rho}
p_n(e)\approx 1-\text{pchisq}\left(\frac{\left(\frac{2}{n}\sum_{i=1}^n
\rho^{(2)}\left(\frac{r_i(e)}{\sigma_n}\right)\right)(s_{\rho,n}(e)-s_{\rho,n}
(e_{\text f}))}{\frac{1}{n}\sum_{i=1}^n\rho^{(1)}
\left(\frac{r_i(e)}{\sigma_n}\right)^2},k-k(e)\right)\,.
\end{equation}
 As the tuning constant $c$ tends to zero the terms
 $\rho_c^{(1)}\left(\frac{r_i(e)}{\sigma_n}\right)^2$ and
 $\rho_c^{(2)}\left(\frac{r_i(e)}{\sigma_n} \right)$ become 
one and zero respectively and the approximation breaks down. 

The results for the stack loss data with $c=1.5$ are given in
Table~\ref{tab:stackloss_M_1.5} and may be compared with those given
in Table~\ref{tab:stackloss_L_1} for the $L_1$-functional.
\begin{table}[h]
\begin{tabular}{ccccccccc}
functional&0&1&2&3&4&5&6&7\\
$P$-value&0.000&0.012&0.000&0.293&0.000&0.006&0.000&1.000\\
$P$-value&3.19e-7&1.23e-2&4.40e-4&3.03e-1&3.79e-8&5.74e-3&3.59e-5&1.00\\
$P$-value&1.89e-6&9.96e-3&1.81e-3&2.33e-1&4.63e-11&2.67e-3&6.66e-5&1.00\\
\quad\\
\end{tabular}
\caption{Encoded $M$-functionals (c=1.5)  and $P$-values for the stack loss
  data based on 5000 simulations: first row the raw values, second row
  the values based on the $\Gamma$-approximation
  (\ref{equ:p_gamma_approx_1}), third row the values based on the
  asymptotic approximation (\ref{equ:asymp_p_rho}). \label{tab:stackloss_M_1.5}} 
\end{table}

\subsection{Least squares regression}
The $L_2$-regression functionals are a special case of the
$M$-functionals for a sufficiently large tuning constant $c$. 
The $P$-values can either be estimated directly using simulations or
using the $\Gamma$-approximation  (\ref{equ:p_gamma_approx_1}) or
using the asymptotic approximation (\ref{equ:asymp_p_rho}) which takes the
form
\begin{equation} \label{equ:lsq_functional_4}
p_n(e)\approx 1-\text{pchisq}\left(\frac{n(\Vert
    {\mathbold y}_n-{\mathbold x}_n(e){\mathbold \beta}_{2,n}(e)\Vert_2^2-\Vert {\mathbold
      y}_n- {\mathbold x}_n{\mathbold \beta}_{2,n}\Vert_2^2)}{\Vert 
  {\mathbold y}_n-{\mathbold x}_n(e){\mathbold \beta}_{2,n}(e)\Vert_2^2},k-k(e)\right)\,.
\end{equation}
If a lower bound $\alpha$ is given for $p_n(e)$ then
(\ref{equ:lsq_functional_4}) is asymptotically
equivalent to the $F$-test in the linear regression model for
testing the null hypothesis that the coefficients of the covariates
not included are zero. 

The $P$-values for the stack loss data are given in
Table~\ref{tab:stackloss_L_2}. 
\begin{table}[h]
\begin{tabular}{ccccccccc}
functional&0&1&2&3&4&5&6&7\\
$P$-value&0.000&0.011&0.000&0.427&0.000&0.005&0.000&1.000\\
$P$-value&1.89e-4&1.07e-2&1.28e-4&4.35e-1&2.80e-5&4.42e-3&4.24e-4&1.00\\
$P$-value&2.53e-4&1.51e-2&1.04e-4&3.11e-1&8.14e-5&4.69e-3&2.20e-4&1.00\\
\quad\\
\end{tabular}
\caption{Encoded $L_2$-functionals and $P$-values for the stack loss
  data based on 5000 simulations: first row the raw values, second row
  the values based on the $\Gamma$-approximation
  (\ref{equ:p_gamma_approx_1}), third row the values based on the
  asymptotic approximation (\ref{equ:lsq_functional_4}).\label{tab:stackloss_L_2}}
\end{table}

 \subsection{Non-linear regression}
 The ideas can be applied mutatis mutandis to non-linear
 regression 
\begin{equation} \label{equ:non_lin}
T_{\text{nl},1,e}(\ep_n)=\argmin_{{\mathbold \beta}(e)} \Vert {\mathbold y}_n-g({\mathbold
  x}_n(e),{\mathbold \beta}(e))\Vert_1
\end{equation} 
with corresponding definitions  for $T_{\text{nl},\rho,e}$ and
$T_{\text{nl},2,e}$. The computational cost is much higher so that
only small values of $k$ are possible.

\subsection{Lower bounds for $P$-values} \label{sec:eval_p}
In contrast to AIC and BIC the $P$-values do not order the different
functionals. One possibility it to choose a cut-off value $p_0(n,k)$ for
$p_n(e)$ and consider only those functionals $T_e$ with $p_n(e)\ge p_0(n,k)$.  A
possible value for $p_0(n,k)$ can be obtained by considering the size of
the $P$-values when all covariates are noise, ${\mathbold
x}_n={\tilde {\mathbold Z}_n}$.  For each such ${\mathbold
   x}_n={\tilde {\mathbold Z}_n}$ the minimum value of $p_n(e)$ over all $e$ can be
 calculated and then simulated for different ${\tilde {\mathbold
   Z}_n}$. The $\alpha$ quantile with for example $\alpha=0.05$  can
then be taken as the value of $p_0(n,k)=p_0(n,k,\alpha)$.  

The minimum of the $p(e)$-values can only be determined by simulation
and then further simulations are required in order to determine the
quantiles of the minimum values. For $L_1$-functionals for the stack
loss data with $n=21$ and $k=3$ the time required with 1000
simulations for each $p(e)$ and 2000 simulations for the quantile was
10 minutes using the approximation to $\vert x\vert$ based on the
Huber $\rho$-function (\ref{equ:huber_rho}) with tuning constant
$c=0.01$. The results using the $\Gamma$-approximation
(\ref{equ:p_gamma_approx_1}) are given in the first line of
Table~\ref{tab:p0nkalpha}.

The computational time for the low birth weight data with $k=9$ was
considerably higher. The time required for 1000 simulations for the
minimum values of the $p_n(e)$ each of which was also based on 500
simulations was 34 hours. The results using
the $\Gamma$-approximation (\ref{equ:p_gamma_approx_1}) are given in
the second line of Table~\ref{tab:p0nkalpha}.

The corresponding $p_0$ values for the $M$-functional and the
$L_2$-functional are based on 500 simulations with 250 simulations for
each $p_n(e)$ value. The computing time for the low birth weight data was
2 1/2 hours in each case. The results are given in
the lines 3-6 of Table~\ref{tab:p0nkalpha}. 
\vspace*{-0.5cm}
\begin{center}
\begin{table}
\begin{tabular}{lll}
$L_1$-functional&&\\
$p_0(21,3,0.01)=0.00368$&$p_0(21,3,0.05)=0.0155$&$p_0(21,3,0.10)=0.0340$\\
$p_0(189,9,0.01)=0.00044$&$p_0(189,9,0.05)=0.0031$&$p_0(189,9,0.10)=0.0068$\\
$L_2$-functional&&\\
$p_0(21,3,0.01)=0.0055$&$p_0(21,3,0.05)=0.0193$&$p_0(21,3,0.10)=0.036$\\
$p_0(189,9,0.01)=0.00011$&$p_0(189,9,0.05)=0.0020$&$p_0(189,9,0.10)=0.0056$
\end{tabular}
\caption{The $P$-values $p_0(n,k,\alpha)$ for the $L_1$- and
  $L_2$-functionals.\label{tab:p0nkalpha}}
\end{table}
\end{center}

For the $M$-functionals with a not too smaller a value of $c$ in
(\ref{equ:huber_rho}) and the $L_2$-functionals use can be made of the
approximations (\ref{equ:asymp_p_rho}) and
(\ref{equ:lsq_functional_4}) respectively which allow simulations 
for larger values of $k$. The computational load can be further reduced as
follows. Without loss of generality suppose $\Vert {\mathbold y}_n
\Vert_2=1$. As all the  ${\mathbold x}_n={\tilde {\mathbold Z}_n}$ are
standard Gaussian white noise random variables it follows for the
$L_2$-functionals
\begin{equation} \label{equ:lsq_functional_2.1}
\frac{\Vert
{\mathbold Z}_n(e^c)^t({\mathbold y}_n-{\mathbold x}_n(e){\mathbold
  \beta}_{2,n}(e))\Vert_2^2}{\Vert  
  {\mathbold y}_n-{\mathbold x}_n(e){\mathbold
    \beta}_{2,n}(e)\Vert_2^2}\approx \Vert {\mathbold
  Z}_n(e^c)^t{\mathbold 
  y}_n\Vert_2^2\approx\sum_j\Vert {\mathbold Z}_n(e^c_j)^t{\mathbold y}_n\Vert_2^2 
\end{equation}
where $e^c=\sum_je^c_j$ and each $e^c_j$ has only one element not
equal to zero. Furthermore the ${\mathbold Z}_n(e^c_j)^t{\mathbold
  y}_n$ are independent $N(0,1)$ random variables. This yields the
asymptotic approximation 
\begin{equation} \label{equ:chipr}
{\tilde p}(e)=1-\text{pchisq}\Big(\sum_{j \in S}\chi^2_1(j),\vert S\vert\Big)
\end{equation}
where the $\chi^2_1(j),j=1,\ldots,k$ are independent $\chi^2$ random
variables with one degree of freedom.  Taking the minimum over $e$, simulating
sets of $\chi^2$ random variables and then taking the $\alpha$
quantile gives a value ${\tilde p}_0(k,\alpha)$. It is only necessary
to perform the simulations one for each value of $k$. The ${\tilde
  p}_0(k,\alpha)$ can be approximated by
\begin{equation} \label{equ:min_pr_approx}
{\tilde p}_0(k,\alpha)\approx \exp(c_1(k)+
c_2(k)\log(\alpha)+c_3(k)\log(\alpha)^2)
\end{equation}
for $\alpha<0.5$ (see Chapter 2.9 of \cite{DAV14}) and can be used in
place of the $p_0(n,k,\alpha)$. For $k=17$ the time required on a
standard laptop was 7 hours 42 minutes. As an example
\[{\tilde p}_0(9,0.01)=0.00028,\,{\tilde p}_0(9,0.05)=0.0025,\, {\tilde
 p}_0(9,0.10)=0.0059\,.\]
The results compare well with those based on simulations as given in
Table~\ref{tab:p0nkalpha} and their computing costs are essentially
zero. This suggests that they can be used as guidelines when
simulations are too expensive.

\subsection{Choosing functionals}  \label{sec:choose_func}
In view of the interpretation of the $P$-value of a functional the
first step is to decide on a cut-off value $p_0$ and then restrict 
consideration to those functionals $T_e$ with $p_n(e)>p_0$. A possible choice
of $p_0$ is $p_0=p_0(n,k,\alpha)$ and this will be done below. The
choice may be further restricted by requiring that for each such $e$
and for all $e'$ with $e'<e$ pointwise all the omitted covariates are
relevant with respect to $e$. More precisely $p_n(e',e) < p_0(n,k(e),\alpha)$ for
all $e'< e$ where $p_n(e',e)$ is the $P$-value of $e'$ calculated with
respect to the covariates ${\mathbold x}(e)$.  A final choice may be
made by choosing that functional $T_e$ with the highest $p_n(e)$-value.

The results of applying the above strategy to the stack loss data with
the $L_1$-norm are as follows. Taking the cut-off value to be
$p_0(21,3,0.01)=0.00368$ (based on Table~\ref{tab:p0nkalpha}) the first step
results in the functionals based on $e_1=(1,0,0)$, $e_2=(1,1,0)$ and
$e_3=(1,0,1)$. The second step eliminates the functionals based on
$e_2$ and $e_3$. The choices $\alpha=0.05$ and $\alpha=0.1$ both lead
to $e=e_2$ with $p_n(e_2)=0.232$.

The results for the low birth weight data are as follows. The first
step with $\alpha=0.01$ leads to 379 functionals. The second step
results in the single functional encoded as 260 with just two
variables and a $P$-value of $4.25e$-3. Putting 
$\alpha=0.05$ results in 221 functionals after the first step. The
second step yields the five functionals encoded as 292, 166,
260, 36 and 60. The functional 292 has the highest $P$-value with
$p_n(e)=0.074$. Finally the choice $\alpha=0.1$ results in 149
functionals after the first step. The second step reduces this to the
seven functionals 308, 292, 166, 262, 38, 52 and 260 of which the
functional encoded as 308 has the highest $P$-value equal to 0.321.

The strategy described above `guarantees' that all included covariates
are relevant. If it is more important not to exclude covariates which may
have an influence at the possible cost of including some irrelevant
covariates the this may be done by increasing the value of $\alpha$ in
$p_0(n,k,\alpha)$ or by simply specifying some cut-off level $p_0$
judged to be appropriate.

Although AIC and BIC list the models in order of preference they give
no indication as to whether any of the models under consideration is
an adequate approximation to the data or not. Presumably this is the
responsibility of the user before  applying the criterion.
The first ten models for the birth weight data based on BIC are
encoded as 
\begin{equation} \label{equ:encoded_vals}
308, 310, 436, 438, 294, 292, 182, 316, 422,309\,.
\end{equation}
The functionals obtained using the $P$-values strategy are encoded as
260, 292 and 308. Their positions in the BIC list are 66, 6 and 1
respectively. 
The second model on the BIC list is encoded as 310 and includes the
additional covariate `weight of mother' compared to the 308
model. If one uses  \cite{R13} to do an $L_1$  regression based on
the covariates corresponding to 310 the 95\%  confidence interval for
`weight of mother' includes zero. This may be interpreted as a
non-significant effect given the other covariates. This interpretation
is consistent with the $P$-value strategy with $\alpha=0.1$ where the
encoded value 310 is not included in the list
(\ref{equ:encoded_vals}). The reason is that the $P$-value 
for the functional excluding `weight of mother' has a $P$-value of
0.112 which exceeds $p_0(189,5,0.1)=0.017$. 

This shows that models may be high in the BIC list although they
contain variables which are not significantly better than random
noise. This can be made more explicit by replacing the all covariates
by random noise and using simulations to determine how often a model
containing a random noise covariate is first on the BIC list.  This
was simulated 500 times with the weight of the child as the dependent
variable. This happened in 43\% of the cases. With $\alpha=0.1$ the
$P$-value strategy is calibrated to do this in $100\alpha\%=10\%$ of
the cases. The simulations resulted in 9\%.

\section{Non-Significance regions}
\subsection{The median and $M$-functionals}
The 0.95-non-significance region for the median of the stack loss data was
defined and calculated in Section~\ref{sec:n-s-r_example} with the
result $[11.86, 18.71]$.  In general the $\alpha$-non-significance region
is defined by
${\mathbold y}_n$ is  
\begin{eqnarray} 
\lefteqn{\mathcal{NS}({\mathbold y}_n,\text{med},\alpha)}\label{equ:med_approx_gen}\\
&=&\left\{m:\sum_{i=1}^{n}\vert
y_i-m\vert-\sum_{i=1}^{n}\vert y_i-\text{med}({\mathbold y}_n)\vert\le
\text{ql1}(\alpha,m,{\mathbold y}_n)\right\}\nonumber 
\end{eqnarray}
where  $\text{ql1}(\alpha,m,{\mathbold y}_n)$ is the $\alpha$-quantile
of
\begin{equation} \label{equ:med_approx0_gen}
\sum_{i=1}^n\vert y_i-m\vert-\inf_b\,\sum_{i=1}^n\vert
y_i-m-bZ_i\vert
\end{equation}
and the $Z_i$ are standard Gaussian white noise.

The non-significance region (\ref{equ:med_approx_gen}) can be calculated
as follows. Put
\begin{equation} \label{equ:quad_approx1}
f(m,\alpha,{\mathbold y}_n)=\text{ql1}(\alpha,m,{\mathbold y}_n)-
\sum_{i=1}^n\vert y_i-m\vert+\sum_{i=1}^n\vert y_i -
\text{med}({\mathbold y}_n)\vert
\end{equation}
and note that $f(\text{med}({\mathbold y}_n),\alpha,{\mathbold y}_n)\ge 0$.
Now determine an order statistic $y_{(nl)}$ with
$nl=\text{qbinom}((1-\beta)/2,n,0.5))$ for a suitably large $\beta$
such that $f(y_{(nl)},\alpha,{\mathbold y}_n)< 0$. Interval bisection
combined with simulations can now be used to find an approximate
solution $m_{\text{lb}}$ of  $f(m,\alpha,{\mathbold y}_n)
=0$. This gives a lower bound and the same process can be used to get
an upper bound $m_{\text{ub}}$ to give $\mathcal{NS}({\mathbold
  y}_n,\text{med},\alpha)=[m_{\text{lb}},m_{\text{ub}}]$.

Non-significance regions for $M$-functionals $T_{\rho}$ are defined
analogously by replacing (\ref{equ:med_approx_gen}) by
\begin{eqnarray} 
\lefteqn{\mathcal{NS}({\mathbold y}_n,T_{\rho},\alpha)}\label{equ:m_approx_gen}\\
&=&\left\{m:\sum_{i=1}^{n}\rho\left(\frac{y_i-m}{\sigma_n}\right)-
  \sum_{i=1}^{n}\rho\left(\frac{y_i-T_{\rho}(\ep_n)}{\sigma_n}\right)\le
\text{qrho}(\alpha,m,{\mathbold y}_n)\right\}\nonumber 
\end{eqnarray}
where  $\text{qrho}(\alpha,m,{\mathbold y}_n)$ is the $\alpha$-quantile
of
\begin{equation} \label{equ:m_approx0_gen}
\sum_{i=1}^n\rho\left(\frac{y_i-m}{\sigma_n}\right)-
\inf_b\,\sum_{i=1}^n\rho\left(\frac{y_i-m-bZ_i}{\sigma_n}\right),
\end{equation}
the $Z_i$ are standard Gaussian white noise and $\sigma_n$ is
a scale functional whose default value in this situation is
(\ref{equ:default_sig_2}).

For smooth functions $\rho$ an asymptotic expression for the
non-significance region is available. Let $\rho^{(1)}$ and
$\rho^{(2)}$ denote the first and second derivative of $\rho$. A
Taylor series expansion results in 
\begin{equation}\label{equ:m_approx_gen_asymp_1}
\mathcal{NS}({\mathbold y}_n,T_{\rho},\alpha)\approx \left\{m:\vert
T_{\rho}(\ep_n)-m\vert \le
\text{qnorm}((1+\alpha)/2)\sigma_n\sqrt{v(T_{\rho},\ep_n)/n}\,\right\}
\end{equation}
where
\begin{equation}\label{equ:m_approx_gen_asymp_2}
v(T_{\rho},\ep_n)= \frac{\frac{1}{n}\sum_{i=1}^n\rho^{(1)}
  \left(\frac{y_i-T_{\rho}(\ep_n)}{\sigma_n}\right)^2}
{\left(\frac{1}{n}\sum_{i=1}^n\rho^{(2)}\left(\frac{y_i-T_{\rho}(\ep_n)}
{\sigma_n}\right)\right)^2}\,.
\end{equation}
This latter expression is well known in robust statistics and
corresponds to the asymptotic variance of an $M$-location
functional: the non-significance region
(\ref{equ:m_approx_gen_asymp_1}) is the corresponding
$\alpha$-confidence region for the `unknown' $T_{\rho}(P)$. In the
special case $\rho(u)=u^2/2$ (\ref{equ:m_approx_gen_asymp_1}) is 
the asymptotic $\alpha$-confidence region for the mean based on
Gaussian errors.

\subsection{$L_1$ regression} \label{sec:approx_qunatile}
 The idea carries over to the $L_1$ regression functional. For any
 ${\mathbold \beta}$ put
\begin{equation}  \label{equ:rq_approx_1}
\Gamma({\mathbold y}_n,{\mathbold x}_n,{\mathbold \beta},{\mathbold Z}_n)=\Vert {\mathbold
  y}_n-{\mathbold x}_n{\mathbold \beta}\Vert_1-\inf_{\mathbold b}\Vert 
{\mathbold y}_n-{\mathbold x}_n{\mathbold \beta}-{\mathbold
  Z}_n{\mathbold b}\Vert_1
\end{equation}
and denote the $\alpha$-quantile of $\Gamma({\mathbold y}_n,{\mathbold
  x},{\mathbold \beta},{\mathbold Z}_n)$ by $\text{q1}(\alpha,{\mathbold \beta},{\mathbold
  y}_n,{\mathbold x}_n)$. An $\alpha$-non-significance region is then
  defined as
\begin{equation} \label{equ:rq_approx_2}
\mathcal{NS}({\mathbold y}_n,{\mathbold x}_n,\alpha,T_{1})=\{{\mathbold \beta}:\Vert
{\mathbold y}_n-{\mathbold x}_n{\mathbold \beta}\Vert_1-\Vert {\mathbold y}_n-{\mathbold
  x}_n{\mathbold \beta}_{1,n}\Vert_1 \le \text{q1}(\alpha,{\mathbold \beta},{\mathbold
  y}_n,{\mathbold x}_n)\}
\end{equation}
where ${\mathbold \beta}_{1,n}=T_{1}(\ep_n)$.

As it stands the non-significance region is difficult to calculate as it
requires a grid of values for the possible values of ${\mathbold \beta}$ and the
values of  $\text{q1}(\alpha,{\mathbold \beta},{\mathbold y}_n,{\mathbold x}_n)$
have to be estimated using simulations. If the quantiles are largely
independent of the ${\mathbold \beta}$-values then
$\text{q1}(\alpha,{\mathbold \beta}_{1},{\mathbold  y}_n,{\mathbold
  x}_n)$ can be used with a large reduction in computation.
Section~\ref{sec:cov_prop} contains some asymptotics which suggest that
the independence may hold for large  sample sizes $n$. The defining
inequality in (\ref{equ:rq_approx_2}) will still have to be checked
over a grid of values.

Most software packages provide only confidence regions for the
individual components of ${\mathbold \beta}$. Corresponding component wise
non-significance regions can be defined with a large
reduction in the computational overload. For the first component
$\beta_1$ of $T_{1}(\ep_n)$ the $\alpha$-non-significance region is
given by 
\begin{eqnarray}
\lefteqn{\mathcal{NS}({\mathbold y}_n,{\mathbold x}_n,\alpha,
  T_{1,1})=\Big\{\beta_1:\inf_{\beta_2,\ldots,\beta_k}\Big\Vert
  {\mathbold y}_n-{\mathbold x}_{\cdot 
  1}\beta_1-\sum_{j=2}^k{\mathbold x}_{\cdot j}\beta_j\Big\Vert_1} \label{equ:rq_approx_3}\\
&&\hspace{3cm}-\Big\Vert
{\mathbold y}_n-{\mathbold x}_n{\mathbold \beta}_{1,n}\Big\Vert_1 \le
\text{q1}(\alpha,\beta_1,{\mathbold
  y}_n,{\mathbold x}_n)\Big\}\nonumber
\end{eqnarray}
where $\text{q1}(\alpha,\beta_1,{\mathbold y}_n,{\mathbold x}_n)$ is
the $\alpha$-quantile of
\begin{eqnarray} 
\Gamma_1({\mathbold y}_n,{\mathbold x}_n,\beta_1,Z_{\cdot 1})
&=&\inf_{\beta_2,\ldots,\beta_k}\Big\Vert {\mathbold y}_n-{\mathbold x}_{\cdot
  1}\beta_1
-\sum_{j=2}^k{\mathbold x}_{\cdot j}\beta_j\Big\Vert_1-\label{equ:rq_approx_4}\\
&&\inf_{b_1,\beta_2,\ldots,\beta_k}\Big\Vert
{\mathbold y}_n-{\mathbold x}_{\cdot 
 1}\beta_1-\sum_{j=2}^k{\mathbold x}_{\cdot j}\beta_j-{\mathbold
 Z}_{\cdot 1}b_1\Big\Vert_1 \nonumber 
\end{eqnarray}

The non-significance intervals of the stack loss data and for
comparison the 0.95-confidence intervals are given in
Table~\ref{tab:approx_stackloss2}.
\begin{table}[h]
\begin{tabular}{cccc}
&Air.Flow&Water.Temp&Acid.Conc\\
Non-sig. intervals (\ref{equ:rq_approx_3})&(0.552,1.082)&(0.225,1.603)&(-0.345,0.102)\\
$L_1$ confidence intervals&(0.509,1.168)&(0.272,3.037)&(-0.278,0.015)\\
\quad\\
\end{tabular}
\caption{First line: 0.95-non-significance intervals for the stack loss
  data. Second line: 0.95-confidence intervals produced by \cite{KOEN10} for the
default choice `se=rank'. \label{tab:approx_stackloss2}}
\end{table}

\subsection{$M$-regression functionals}
Non-significance regions for $M$-regression functionals are defined in
the same manner as for $L_1$ regression. Just as in
Section~\ref{sec:M_modchc} the computational burden can be reduced for
large $n$ by using the asymptotic expressions. These result in 
\begin{eqnarray} 
\lefteqn{\hspace{-1cm}\mathcal{NS}({\mathbold y}_n,{\mathbold
    x}_n,\alpha,T_{\rho})=\Biggl\{{\mathbold
    \beta}:\sum_{i=1}^n\rho\left(\frac{y_i-{\mathbold
        x}_{i\cdot}^t{\mathbold \beta}}{\sigma_n}\right)-\sum_{i=1}^n
  \rho\left(\frac{y_i-{\mathbold x}_{i\cdot}^t{\mathbold
        \beta_{\rho}}}{\sigma_n}\right)\le }\label{equ:m_approx_1}\\ 
&&\hspace{3cm}\frac{\text{qchisq}(\alpha,k)}{2}\frac{\sum_{i=1}^n
\rho^{{(1)}^2}\left(\frac{y_i-{\mathbold x}_{i\cdot}^t{\mathbold \beta}}{\sigma_n}\right)}{\sum_{i=1}^n
\rho^{(2)}\left(\frac{y_i-{\mathbold x}_{\cdot}i^t{\mathbold \beta}}{\sigma_n}\right)}\Biggr\}\nonumber
\end{eqnarray}
where ${\mathbold \beta_{\rho}}=T_{\rho}(\ep_n)$. This can be further
simplified to
\begin{eqnarray}
\lefteqn{\hspace{1.5cm}\mathcal{NS}({\mathbold y}_n,{\mathbold
    x}_n,\alpha,T_{\rho})=\Bigg\{{\mathbold \beta}:({\mathbold
    \beta}-{\mathbold \beta}_{\rho})^t{\mathbold x}_n^t{\mathbold
    x}_n({\mathbold \beta}-{\mathbold \beta}_{\rho}) 
  \le } \label{equ:m_approx_2}\\
&&\hspace{2cm}\hspace{3cm}\left.\text{qchisq}(\alpha,k)\frac{\sum_{i=1}^n
\rho^{{(1)}^2}\left(\frac{y_i-{\mathbold x}_{i\cdot}^t{\mathbold
      \beta_{\rho}}}{\sigma_n}\right)}{\sum_{i=1}^n 
\rho^{(2)}\left(\frac{y_i-{\mathbold x}_{\cdot}i^t{\mathbold
      \beta_{\rho}}}{\sigma_n}\right)}\right\}\,.\nonumber 
\end{eqnarray}

\subsection{Least squares regression}
The method goes through for the least squares functional with the
advantage that explicit expressions are available. The result
corresponding to (\ref{equ:m_approx_1}) is
\begin{eqnarray} 
\lefteqn{\hspace{-2cm}\mathcal{NS}({\mathbold y}_n,{\mathbold x}_n,\alpha,T_2)=\Bigg\{{\mathbold \beta}:\Vert {\mathbold y}_n-{\mathbold x}_n{\mathbold \beta}\Vert_2^2-\Vert
{\mathbold y}_n-{\mathbold x}_n{\mathbold \beta}_{2,n}\Vert_2^2\le} \label{equ:lsq_approx_1}\\
&&\hspace{2cm}\frac{\Vert {\mathbold y}_n-{\mathbold
     x}_n{\mathbold \beta}\Vert_2^2}{n}\text{qchisq}(\alpha,k)\Bigg\}
\nonumber
\end{eqnarray}
which is the same as
\begin{eqnarray}
\lefteqn{\hspace{-1.5cm}\mathcal{NS}({\mathbold y}_n,{\mathbold
    x}_n,\alpha,T_{2})=\Bigg\{{\mathbold \beta}:({\mathbold
    \beta}-{\mathbold \beta}_{2,n})^t{\mathbold x}_n^t{\mathbold
    x}_n({\mathbold \beta}-{\mathbold \beta}_{2,n}) 
  \le } \label{equ:lsq_approx_2}\\
&&\hspace{2cm}\frac{\Vert
{\mathbold y}_n-{\mathbold x}_n{\mathbold \beta}_{2,n}\Vert_2^2 
\text{qchisq}(\alpha,k)}{n-\text{qchisq}(\alpha,k)}\Bigg\}\,. \nonumber
\end{eqnarray}
where ${\mathbold \beta}_{2,n}=T_2(\ep_n)$. The region is
asymptotically equivalent to a standard $\alpha$-confidence region for
the `true' parameter value.

\subsection{Covering properties} \label{sec:cov_prop}
The concept of a non-significance region makes no mention of a model
or true values. Nevertheless there are situations where a model and
its parameters are well founded and relate to well-defined properties 
of the real world. In such cases there is an interest in
specifying a region which includes the real world value with the
required frequency in repeated measurements. It has to
be kept in mind however that covering true parameter values in
simulations is not the same as covering the corresponding real values
for real data (see Chapter 5.5 of \cite{DAV14}, \cite{STIG77}, Chapter
8.1 of \cite{HAMRONROUSTA86}, \cite{KUNBERHAMP93}).

Given this there is an interest in the covering properties of
non-significance regions. Table~\ref{tab:l1_approx_reg1} gives the
frequencies with which the non-significance intervals
(\ref{equ:med_approx_gen}) and the confidence intervals based on the rank
statistics cover the population median and 
also the average lengths of the intervals. The results are for the
normal, Cauchy, $\chi^2_1$ and the Poisson ${\mathfrak
  P}{\mathfrak o}(4)$ distributions and four different sample sizes
$n=10, 20, 50, 100$ and are  based on 1000 simulations. The
discreteness of Poisson distribution was taken into account in the
calculations of the non-significance region as follows.  If an
non-significance interval 
$[\ell,u]$ contains an integer it is by $[\lceil
\ell \rceil,\lfloor u\rfloor]$. If it does not contain an integer it is
replaced by $[\lfloor\ell\rfloor,\lceil u\rceil]$. The covering
frequencies and lengths refer to this modified interval.
\begin{table}
\begin{tabular}{ccccccc}
&$n$&10&20&50&100\\
$N(0,1)$& (\ref{equ:med_approx_gen})&0.940 1.512&0.954 1.040&0.948 0.648&0.942 0.464\\ 
&rank&0.968 2.046&0.968 1.198&0.970 0.767&0.964 0.530\\
$C(0,1)$&(\ref{equ:med_approx_gen})&0.960 3.318&0.956 1.670&0.960 0.958&0.952 0.629\\
&rank&0.978 5.791&0.950 1.850&0.968 1.069&0.964 0.700\\
$\chi^2_1$& (\ref{equ:med_approx_gen})&0.944 1.368&0.936 0.877&0.932 0.550&0.942 0.396\\
&rank&0.982 2.064&0.958 1.086&0.970 0.675&0.968 0.452\\
Pois(4)&(\ref{equ:med_approx_gen})&0.934 1.918&0.925 0.993&0.926 0.288&0.938 0.071\\
&rank&0.996 3.948&0.964 2.342&0.997 1.573&1.000 1.085\\
\\
\end{tabular}
\caption{Covering frequencies and average interval lengths based on
  1000 simulations for the   median for the   0.95-non-significance
  intervals as defined by (\ref{equ:med_approx_gen}) and
  (\ref{equ:med_approx0_gen}) with $Z=N(0,1)$ and the 0.95-confidence
  intervals based on the ranks. For each sample size the first column
  gives the covering frequency and the second the average interval
  length. \label{tab:l1_approx_reg1}}    
\end{table}
In this well defined situation Table~\ref{tab:l1_approx_reg1}
indicates that the 0.95-non-significance intervals also have covering
probabilities of about 0.95. The finite sample behaviour seems to be
better than that of the ranks procedure. Both methods have
approximately the correct covering frequencies but the lengths of the
non-significance intervals are uniformly smaller than the lengths of the
confidence intervals. 

There is some theoretical explanation as to why the non-significance
regions have covering frequencies given by $\alpha$, at least
asymptotically. Consider firstly i.i.d. integer valued random
variables $Y_j$ with a unique median $\nu$. Then for a large sample
size $n$
\[ \sum_{j=1}^n \vert Y_j-\nu-bZ_j\vert\] 
is, with large probability, minimized by putting $b=0$. In other words
the 0.95-non-significance interval is simply $[\nu,\nu]$ with a covering
probability tending to one.  This is illustrated by the Poisson
distribution in Table~\ref{tab:l1_approx_reg1}.

Suppose that the $Y_j$ are continuous random variables with median 0
and a density $f$ which is continuous at 0 with $f(0)>0$. Then the
approximation
\begin{equation} \label{equ:rq_cov_prob_1}
\sum_{i=1}^n \left\vert Y_i-\frac{bZ_i}{\sqrt{n}}\right\vert \approx
\sum_{i=1}^n\vert Y_i\vert -bN(0,1)+f(0)b^2
\end{equation}
holds (see the Appendix for a heuristic proof) and minimizing over $b$ gives 
\begin{equation}  \label{equ:rq_cov_prob_2}
\inf_b\,\sum_{i=1}^n \left\vert Y_i-\frac{bZ_i}{\sqrt{n}}\right\vert \approx
\sum_{i=1}^n\vert Y_i\vert- \frac{\chi^2_1}{4f(0)}\,.
\end{equation}
Moreover the same proof gives
\begin{equation} \label{equ:rq_cov_prob_3}
\inf_b\,\sum_{i=1}^n \left\vert Y_i-\text{med}({\mathbold
    Y}_n)-\frac{\theta}{\sqrt{n}}-\frac{bZ_i}{\sqrt{n}}\right\vert
\approx \sum_{i=1}^n\vert Y_i-\text{med}({\mathbold
    Y}_n)\vert+f(0)\theta^2-\frac{\chi^2_1}{4f(0)}
\end{equation}
from which the asymptotic $\alpha$-non-significance interval 
\begin{equation} \label{equ:rq_cov_prob_4}
\left[\text{med}({\mathbold
  Y}_n)-\sqrt{\frac{\text{qchisq}(\alpha,1)}{4f(0)^2n}},\,
  \text{med}({\mathbold Y}_n)+\sqrt{\frac{\text{qchisq}(\alpha,1)}{4f(0)^2n}}\,\right]
\end{equation}
as defined in (\ref{equ:med_approx_gen}) and
(\ref{equ:med_approx0_gen}) follows. This latter interval is 
the same as the asymptotic confidence interval based on the
median. Just as for the inverse rank method it does not require an
estimate of $f(0)$.

$L_1$ linear regression can be treated in the same
manner. Corresponding to  (\ref{equ:rq_cov_prob_1}) one has
\begin{equation} \label{equ:rq_cov_prob_5}
  \sum_{i=1}^n \left\vert Y_i-\frac{{\mathbold Z}_{i\cdot}^t{\mathbold
b}}{\sqrt{n}}\right\vert \approx 
 \sum_{i=1}^n\vert Y_i\vert -N(0,\mathbold{I}_k)^t{\mathbold b}+f(0)\Vert
 {\mathbold b}\Vert_2^2\,.
\end{equation}
Applying this to the $L_1$ regression functional gives
\begin{equation}  \label{equ:rq_cov_prob_6}
\inf_{\mathbold b}\,\sum_{i=1}^n \left\vert Y_i-{\mathbold
    x}_{i\cdot}^t{\mathbold \beta}_{1,n}-\frac{{\mathbold
      x}_{i\cdot}^t{\mathbold 
      \theta}}{\sqrt{n}}-\frac{{\mathbold Z}_{i\cdot}^t{\mathbold b}}{\sqrt{n}}\right\vert 
\approx\sum_{i=1}^n \vert Y_i-{\mathbold x}_{i\cdot}^t{\mathbold \beta}_{1,n}\vert
+f(0){\mathbold \theta}^t{\mathbold Q}_n{\mathbold \theta}-\frac{\chi^2_k}{4f(0)^2}
\end{equation}
where ${\mathbold Q}_n=\frac{1}{n}{\mathbold x}_n^t{\mathbold x}_n$. From this the asymptotic
$\alpha$-non-significance region
\begin{equation}  \label{equ:rq_cov_prob_7}
({\mathbold \beta}-{\mathbold \beta}_{1,n})^t{\mathbold
  Q}_n({\mathbold \beta}-{\mathbold \beta}_{1,n})\le 
\frac{\text{qchisq}(\alpha,k)}{4f(0)^2n}
\end{equation}
follows. It is the same as the $\alpha$-confidence region based on
the $L_1$ regression estimate ${\mathbold \beta}_{1}$, see for example
\cite{ZHOPOR96}.

Table~\ref{tab:cov_stack_loss} gives the covering frequencies and
average interval lengths for data generated according to
\begin{equation} \label{equ:stackloss_data_gen}
Y=-39.69+0.832\cdot Air.Flow+0.574\cdot Water.Temp-0.061\cdot
Acid.Conc +\varepsilon
\end{equation}
using the $L_1$ coefficients for the stack loss data. The
sample size is $n=21$. The following four distributions for
the error term $\varepsilon$ are used:
$\varepsilon=N(0,1)*\text{Res}$, $\varepsilon=\sigma N(0,1)$
for the normal distribution, $\varepsilon=\sigma L^*$ for the Laplace
distribution and $\varepsilon=\sigma C^*$ for the Cauchy distribution
where $L^*$ and $C^*$ are respectively the Laplace and Cauchy
distributions closest to the $N(0,1)$ distribution, $\text{Res}$ are
the residuals and $\sigma$ the mean absolute deviation of the
residuals of the stack loss data.
\begin{table}[h]
\begin{tabular}{ccccc}
&&$\beta_2$&$\beta_3$&$\beta_4$\\
residuals&(\ref{equ:rq_approx_3})&0.944 0.265&0.982 0.682&0.998 0.248\\
&rank&0.976 0.390&0.970 1.205&0.970 0.273\\
Normal&(\ref{equ:rq_approx_3})&0.954 0.381&0.946 1.042&0.964 0.442\\
&rank&0.974 0.435&0.956 1.208&0.962 0.542\\
Laplace&(\ref{equ:rq_approx_3})&0.953 0.501&0.959 1.375&0.952 0.580\\
&rank&0.966 0.594&0.959 1.697&0.960 0.761\\
Cauchy&(\ref{equ:rq_approx_3})&0.928 1.467&0.942 4.052&0.936 1.731\\
&rank&0.936 1.948&0.946 5.676&0.942 2.984\\
\\
\end{tabular}
\caption{Covering frequencies and average interval lengths for data
  generated according to (\ref{equ:stackloss_data_gen}) with different
  distributions for the error term: $\alpha=0.95$.\label{tab:cov_stack_loss}}
\end{table}

Finally, in the case of non-linear $L_1$ regression the asymptotic
$\alpha$-non-significance is, under suitable regularity conditions, 
given by
\begin{equation}  \label{equ:rq_cov_prob_8}
({\mathbold \beta}-{\mathbold \beta}_{{\text nlr1},n})^t{\mathbold
  Q}_n({\mathbold \beta}-{\mathbold \beta}_{{\text nlr1},n})\le 
\frac{\text{qchisq}(\alpha,k)}{4f(0)^2n}
\end{equation}
where 
\[{\mathbold Q}_n=\frac{1}{n}\sum_{i=1}^n {\mathbold
  \nabla}_i{\mathbold \nabla}_i^t\]
and
\[{\mathbold \nabla}_i =\left( \frac{\partial m({\mathbold x}_{i\cdot},{\mathbold \theta})}{\partial
    \theta_1},\ldots,\frac{\partial m({\mathbold x}_{i\cdot},{\mathbold \theta})}{\partial
    \theta_k}\right)^t\,.\] 

\section{Choice of noise} \label{sec:choice_noise}
It is possible to use random variables other than Gaussian. As an
example the 0.95-non-significance intervals for the 
median of the stack loss data using $N(0,1)$, $\pm 1$, $U(-1,1)$, $\pm
\text{beta}(5,5)$ and the standard Cauchy distribution are $(11.88, 18.63)$,
$(11.86,18.25)$  $(11.83, 18.16)$, $(11.93, 18.21)$,  $(11.83, 18.16)$
and $(11.00,18.56)$ respectively. It is clear that the results depend
on the choice of noise to some extent but that at least in this
example the dependence is weak. Given the advantages of Gaussian
noise are the easily available  asymptotic expressions such
(\ref{equ:asymp_p_rho}) it would seem to be the default choice of
noise. 

Other possibilities are to make the noise dependent on the size of the
covariates as in $W_{ij}=x_{ij}Z_{ij}$  or to randomly permute the
covariates  (see \cite{ANDROB01}, \cite{KLI09}).

\section{Appendix}
Consider 
\begin{equation} \label{equ:sum_med}
\sum_{i=1}^n \left\vert
  \varepsilon_i-\frac{{\mathbold U}_{i\cdot}^t{\mathbold b}}{\sqrt{n}}\right\vert
\end{equation}
where the $\varepsilon_i$ are symmetric, i.i.d. random variables with a continuously
differentiable density at $u=0$ with $f(0)>0$.  The ${\mathbold U}_{i\cdot}$ are $k$
dimensional random variables ${\mathbold U}_{i\cdot}=(U_{i,1},\ldots,U_{i,k})^t$ where the
$U_{ij}$ are symmetric i.i.d. random variable with unit variance. The
sum (\ref{equ:sum_med}) may be decomposed as
\begin{eqnarray*}
\sum_{i=1}^n \left\vert
  \varepsilon_i-\frac{{\mathbold U}_{i\cdot}^t{\mathbold b}}{\sqrt{n}}\right\vert
&=& \sum_{\varepsilon_i\le -\left\vert{\mathbold
      U}_{i\cdot}^t{\mathbold
      b}/\sqrt{n}\right\vert}\left(-\varepsilon_i+ {\mathbold
    U}_{i\cdot}^t{\mathbold b}/\sqrt{n}\right) \\
&&+\sum_{\varepsilon_i\ge \left\vert  
{\mathbold U}_{i\cdot}^t{\mathbold  b}/\sqrt{n}\right\vert}\left(\varepsilon_i-{\mathbold
     U}_{i\cdot}^t{\mathbold b}/\sqrt{n}\right)\\ 
&&+\sum_{\vert \varepsilon_i\vert \le \left\vert
    {\mathbold U}_{i\cdot}^t{\mathbold b}/\sqrt{n}\right\vert}\left\vert \varepsilon_i+
  \frac{{\mathbold U}_{i\cdot}^t{\mathbold b}}{\sqrt{n}}\right\vert\\
&=&\sum_{i=1}^n\vert \varepsilon_i\vert
 +\sum_{i=1}^n\pm\frac{{\mathbold U}_{i\cdot}^t{\mathbold
 b}}{\sqrt{n}}
-\sum_{\vert \varepsilon_i\vert \le \left\vert
{\mathbold U}_{i\cdot}^t{\mathbold
  b}/\sqrt{n}\right\vert}\frac{{\mathbold U}_{i\cdot}^t{\mathbold
b}}{\sqrt{n}}\\ 
&&+\sum_{\vert \varepsilon_i\vert \le \left\vert
    \frac{{\mathbold U}_{i\cdot}^t{\mathbold b}}{\sqrt{n}}\right\vert}\left( \left\vert \varepsilon_i-
  \frac{{\mathbold U}_{i\cdot}^t{\mathbold b}}{\sqrt{n}}\right\vert -\vert \varepsilon_i\vert \right)\,.
\end{eqnarray*}
The random variables 
\[V_i=\frac{{\mathbold U}_{i\cdot}^t{\mathbold b}}{\sqrt{n}}\left\{\vert \varepsilon_i\vert \le \left\vert
    \frac{{\mathbold U}_{i\cdot}^t{\mathbold b}}{\sqrt{n}}\right\vert\right\}\]
are i.i.d with mean zero and variance
\[\frac{1}{n}\ex_U \left(({\mathbold U}_{i\cdot}^t{\mathbold b})^2\left(F\left(\left\vert
      \frac{{\mathbold U}_{i\cdot}^t{\mathbold b}}{\sqrt{n}}\right\vert\right)-
  F\left(-\left\vert\frac{{\mathbold U}_{i\cdot}t^{\mathbold b}}{\sqrt{n}}\right\vert\right)\right)
\right)=o\left(\frac{\Vert {\mathbold b}\Vert_2^2}{n}\right)\,.\]    
This together with the central limit theorem implies
\begin{eqnarray*}
\sum_{i=1}^n \left\vert
  \varepsilon_i-\frac{{\mathbold U}_{i\cdot}^t{\mathbold b}}{\sqrt{n}}\right\vert
&=&\sum_{i=1}^n\vert \varepsilon_i\vert +{\mathbold Z}^t{\mathbold b}+\sum_{\vert
  \varepsilon_i\vert \le \left\vert {\mathbold U}_{i\cdot}^t{\mathbold b}/\sqrt{n}\right\vert}\left( \left\vert \varepsilon_i-\frac{{\mathbold U}_{i\cdot}^t{\mathbold b}}{\sqrt{n}}\right\vert -\vert \varepsilon_i\vert\right)+o\left(\Vert {\mathbold b}\Vert_2^2\right)
\end{eqnarray*}
where ${\mathbold Z}\stackrel{D}{=}N(0,I_k)$. Denote the distribution function of
${\mathbold U}_{i\cdot}^t{\mathbold b}$ by $H$. Then
\[\ex_U\left(\left\vert\varepsilon_i-\frac{{\mathbold U}_{i\cdot}^t{\mathbold b}}{\sqrt{n}}\right\vert\left\{\vert
  \varepsilon_i\vert \le \left\vert\frac{{\mathbold U}_{i\cdot}^t{\mathbold b}}{\sqrt{n}}\right\vert\right\}\right)=
\frac{2}{\sqrt{n}}\int_0^{\infty}w\left\{\vert\varepsilon_i\vert \le
  \frac{w}{\sqrt{n}}\right\}\,dH(w) \]
and taking the expected value with respect to $\varepsilon_i$ gives
\begin{eqnarray*}
\lefteqn{\hspace{-2cm}\ex\left(\left\vert\varepsilon_i-\frac{{\mathbold U}_{i\cdot}^t{\mathbold b}}{\sqrt{n}}\right\vert\left\{\vert
  \varepsilon_i\vert \le
  \left\vert\frac{{\mathbold U}_{i\cdot}^t{\mathbold b}}{\sqrt{n}}\right\vert\right\}\right)}\\
&=&
\frac{2}{\sqrt{n}}\int_0^{\infty}w
\left(F\left(\frac{w}{\sqrt{n}}\right)-
  F\left(-\frac{w}{\sqrt{n}}\right)\right)\,dH(w)\\
&\approx&\frac{4f(0)}{n}\int_0^{\infty}w^2\,dH(w)
=\frac{2f(0)\Vert {\mathbold b}\Vert_2^2}{n}
\end{eqnarray*}
as the $U_{ij}$ are symmetric random variables with variance 1. A similar
calculation gives
\[\ex\left(\vert\varepsilon_i\vert\left\{\vert
  \varepsilon_i\vert \le
  \left\vert\frac{{\mathbold U}_{i\cdot}^t{\mathbold b}}{\sqrt{n}}\right\vert\right\}\right)\approx
\frac{f(0)\Vert {\mathbold b}\Vert_2^2}{n}\,.\] 
Putting this together leads to
\[\sum_{i=1}^n \left\vert
  \varepsilon_i-\frac{{\mathbold U}_{i\cdot}^t{\mathbold b}}{\sqrt{n}}\right\vert \approx
\sum_{i=1}^n\vert \varepsilon_i\vert +Z^t{\mathbold b} +\frac{f(0)\Vert
  {\mathbold b}\Vert_2^2}{n}
\]
and minimizing over ${\mathbold b}$ results in
\[\inf_{\mathbold b}\sum_{i=1}^n \left\vert
  \varepsilon_i-\frac{{\mathbold U}_{i\cdot}^t{\mathbold b}}{\sqrt{n}}\right\vert\approx\sum_{i=1}^n\vert
\varepsilon_i\vert -\frac{\chi_k^2}{4f(0)}\]
where $\chi_k^2$ is a chi-squared random variable with $k$ degrees of freedom.

\end{document}